# A new derivative-free optimization method: Gaussian Crunching Search

B.Wong

Techdesignism Machine Learning Laboratory

Abstract: Optimization methods are essential in solving complex problems across various domains. In this research paper, we introduce a novel optimization method called Gaussian Crunching Search (GCS). Inspired by the behaviour of particles in a Gaussian distribution, GCS aims to efficiently explore the solution space and converge towards the global optimum. We present a comprehensive analysis of GCS, including its working mechanism, and potential applications. Through experimental evaluations and comparisons with existing optimization methods, we highlight the advantages and strengths of GCS. This research paper serves as a valuable resource for researchers, practitioners, and students interested in optimization, providing insights into the development and potential of Gaussian Crunching Search as a new and promising approach.

## Introduction

Derivative-free optimization methods play a crucial role in solving complex problems where the objective function is either non-differentiable or computationally expensive to evaluate. In this research paper, we present a novel derivative-free optimization method called Gaussian Crunching Search (GCS) that we have developed to address these challenges.

Gaussian Crunching Search draws inspiration from the natural phenomenon of Gaussian distributions and leverages its principles to optimize the search process. Unlike traditional optimization methods that rely on gradient information, GCS explores the solution space by mimicking the behaviour of particles in a Gaussian distribution. This unique approach allows GCS to efficiently navigate the search landscape, converging towards the optimal solution without the need for derivative information.

The primary objective of this research is to introduce and evaluate the effectiveness of Gaussian Crunching Search as a new derivative-free optimization method. We will delve into the step-by-step process of GCS, providing a comprehensive understanding of its working mechanism. By analyzing it with existing derivative-free optimization techniques, we aim to highlight the advantages and potential applications of GCS.

This research paper aims to contribute to the field of optimization by presenting a detailed description and analysis of Gaussian Crunching Search. We will discuss the motivation behind the development of GCS, its unique features, and its potential benefits in various problem domains. By understanding the inner workings of GCS, researchers, practitioners, and students can gain insights into its applicability and evaluate its potential for solving real-world optimization problems efficiently and effectively.

In the following sections, we will present a thorough exploration of Gaussian Crunching Search, including its experimental results, and comparisons with existing derivative-free optimization methods. Through this research, we aim to provide a solid foundation for further studies and applications of GCS, fostering innovation in the field of optimization and paving the way for advancements in solving complex problems without the need for derivative information.

## Limitations of Traditional Derivative-Free Methods

Traditional derivative-free optimization methods, such as "Nelder-Mead", "Powell" "CG", "BFGS", "L-BFGS-B", "TNC", "COBYLA", "SLSQP" and "Trust-Constr" have long been relied upon for solving complex optimization problems. However, a significant limitation arises when the difference in the result value of the objective function becomes extremely small, exceeding the accuracy that can be effectively expressed by the computer program. This issue poses a significant challenge, as many of these commonly used methods fail to function optimally under such circumstances.

When the difference in the objective function's result value becomes minuscule, the traditional derivative-free methods struggle to accurately converge towards the optimal solution. This limitation is primarily attributed to the numerical precision of the computer program, which can only represent a certain level of accuracy. Consequently, the inability of these methods to handle extremely small differences in the objective function hinders their effectiveness in finding the true optimal solution.

The problem with the traditional derivative-free methods becomes particularly evident when dealing with highly sensitive optimization problems or when the objective function exhibits a plateau-like behaviour. In such cases, even a small difference in the objective function's value can significantly impact the quality of the solution. As a result, the inability of these methods to handle such scenarios limits their applicability and effectiveness in real-world optimization problems where the objective function's value differences can be extremely small.

**EXAMPLE**

$$f(x, y) = -\lambda e^{-\mu \sqrt{x^2 + y^2}} + \lambda \qquad (1)$$

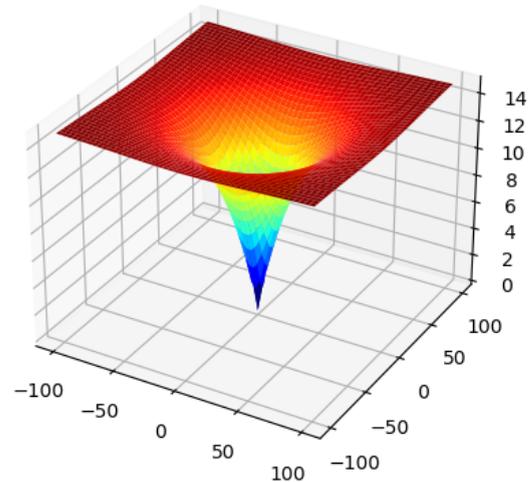

We will evaluate the limitations of traditional commonly used derivative-free methods[1], such as "Nelder-Mead", "Powell", "CG", "BFGS", "L-BFGS-B", "TNC", "COBYLA", "SLSQP", and "Trust-Constr" when faced with small differences in the result value of the objective function that exceeds the accuracy expressible by the computer program. To demonstrate these limitations, we will implement the optimization process using Python and the Scipy library[2]. Specifically, we will consider Eq. (1) as our objective function, which serves as a standard benchmark for testing optimization methods. By executing the optimization process with Python and Scipy, we aim to present a comprehensive analysis of the limitations encountered by these traditional derivative-free methods in handling objective function differences beyond the computational accuracy. Through this analysis, we intend to shed light on the challenges faced by these methods and potentially identify areas for improvement or alternative approaches to address this limitation.

Starting point: [x: 200, y: 200]

| Method | Results | <0.5 |
|---|---|---|
| Nelder-Mead | 2.22E-05 | TRUE |
| Powell | 3.26E-12 | TRUE |
| CG | 15 | FALSE |
| BFGS | 15 | FALSE |
| L-BFGS-B | 15 | FALSE |
| TNC | 2.83E-09 | TRUE |
| COBYLA | 6.42E-05 | TRUE |
| SLSQP | 15 | FALSE |
| trust-constr | 3.42E-09 | TRUE |

Starting point: [x: 400, y: 400]

| Method | Results | <0.5 |
|---|---|---|
| Nelder-Mead | 1.72E-05 | TRUE |
| Powell | 3.60E-12 | TRUE |
| CG | 15 | FALSE |
| BFGS | 15 | FALSE |
| L-BFGS-B | 15 | FALSE |
| TNC | 15 | FALSE |
| COBYLA | 8.15E-05 | TRUE |
| SLSQP | 15 | FALSE |
| trust-constr | 15 | FALSE |

Starting point: [x: 600, y: 600]

| Method | Results | <0.5 |
|---|---|---|
| Nelder-Mead | 15 | FALSE |
| Powell | 15 | FALSE |
| CG | 15 | FALSE |
| BFGS | 15 | FALSE |
| L-BFGS-B | 15 | FALSE |
| TNC | 15 | FALSE |
| COBYLA | 15 | FALSE |
| SLSQP | 15 | FALSE |
| trust-constr | 15 | FALSE |

Python Code:

```python
from numpy import *
from scipy.optimize import minimize, LinearConstraint
import scipy
import numpy as np
np.set_printoptions(precision=2)

def objective(x_list):
    [x,y]= x_list
    ld = 15
    mu = 0.05
    return -ld* exp(-mu* sqrt(x**2 + y**2))  + ld

methods_list = [
"Nelder-Mead"
,"Powell"
,"CG"
,"BFGS"
,"L-BFGS-B"
,"TNC"
,"COBYLA"
,"SLSQP"
,"trust-constr"
]

for i in [200,400,600]:
    start_point = [i,i]
    print('Starting Point: ', [i,i],'\n')
    for x in methods_list:
        print('Method:',x)
        res = minimize(
            objective,
            start_point,
            method=x
        )
        print(res.x)
        print( np.array([res.fun]) )
        print('='*10)
    print('*'*40)
```

We conducted experiments by setting λ to 15 and μ to 0.05 to investigate the performance of traditional commonly used derivative-free methods. The optimization process was initiated with three different starting points: [200,200], [400,400], and [600,600]. Notably, our findings revealed that the majority of the optimization methods faced significant challenges when attempting to converge towards the optimal solution at the [400,400] starting point. Furthermore, it is noteworthy that all of the methods failed to achieve convergence at the [600,600] starting point. These results highlight the limitations of these traditional derivative-free methods in effectively handling optimization problems with particular starting points, emphasizing the need for alternative approaches or improvements in order to overcome such limitations and ensure robust optimization performance.

## Gaussian Crunching Search (GCS)

Traditional approaches to derivative-free optimization have played a crucial role in addressing complex optimization problems. However, these methods encounter a significant drawback when the difference between objective function values becomes exceedingly small, surpassing the precision achievable by computer programs. In light of this limitation, we introduce Gaussian Crunching Search (GCS), an innovative derivative-free optimization technique inspired by the genetic algorithm[3]. GCS aims to alleviate the challenges associated with conventional optimization methods that heavily rely on gradient information. By incorporating a stochastic mutation process, GCS effectively and efficiently explores the solution space, striving to achieve enhanced outcomes across diverse optimization problems.

GCS operates based on the fundamental concept of emulating natural selection and genetic variation. Instead of relying on updates driven by gradients, GCS capitalizes on a stochastic mutation process to induce diversity within the population of solutions. This inherent stochasticity empowers GCS to escape local optima and venture into unexplored regions of the search space. By synergistically combining the exploration capabilities of genetic algorithms with the optimization objectives of derivative-free techniques, GCS provides a promising avenue for attaining more optimal results across a wide spectrum of optimization problems.

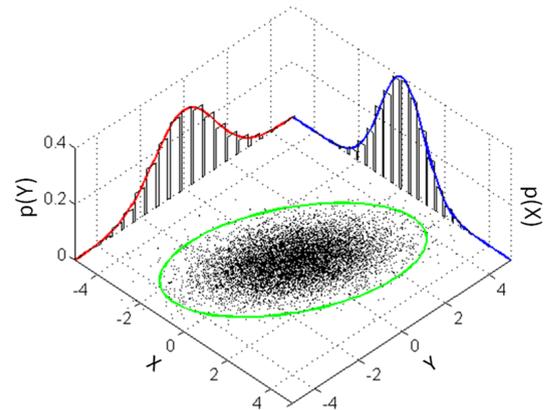

Figure 1: 2D Gaussian Distribution

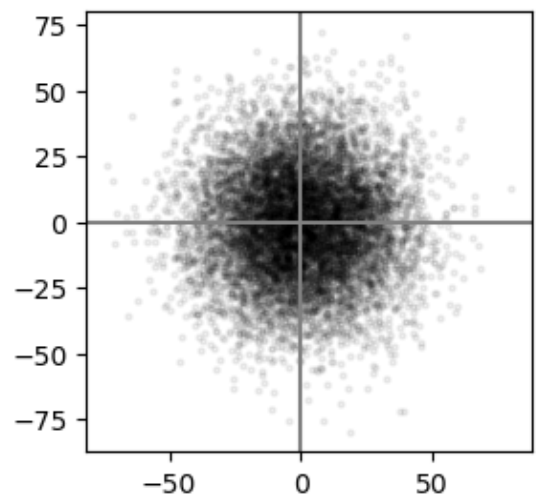

Figure 2: Stochastic Mutation

Gaussian Crunching Search (GCS) optimizes the objective function by employing simulated natural selection. It combines mutation, implemented through a Gaussian distribution, and selection mechanisms in each iteration to enhance the optimization process. The mutation step introduces slight perturbations to the current solution, facilitating the exploration of the search space. Subsequently, the selection mechanism identifies promising solutions based on their fitness values.

A notable aspect of GCS is its dynamic mutation strategy. As iterations progress with the same temporary optimal value, the standard deviation of the mutation distribution gradually increases. This adaptive feature enables GCS to strike a balance between exploration and exploitation. In the early stages of optimization, GCS extensively explores the search space with smaller deviations. As the algorithm approaches a local optimum, the augmented standard deviation encourages further exploration, aiding the discovery of potentially superior solutions in uncharted regions. When a new temporary optimal value is encountered, GCS resets the standard deviation, resulting in variations in the shape of the Gaussian distribution, reminiscent of a "crunching" effect that lends its name to the technique.

Through this iterative process of mutation and selection, GCS effectively harnesses natural selection principles to fine-tune the optimization of the objective function. By dynamically adjusting the mutation strategy based on the number of iterations with the same temporary optimal value, GCS exhibits an inherent ability to adapt the exploration-exploitation trade-off, leading to improved efficiency and effectiveness in optimizing performance.

GCS Pseudocode:

```
SET sd TO 1
SET current_pos TO start_point
SET current_obj_value TO objective_function(current_pos)
SET last_pos TO current_pos

FOR i IN range(10000):

    SET temp_pos TO current_pos + random_number(mean: 0, standard_deviation: sd)

    SET temp_obj_value TO objective(temp_pos)

    IF  temp_obj_value < current_obj_value:
       SET last_pos TO current_pos
       SET current_pos TO temp_pos
       SET current_obj_value TO temp_obj_value

       SET sd TO 1
    ELSE:
       SET sd TO sd * 1.01
       IF sd > the_limitation_of_python:
          SET sd TO 1
```

GCS Python Code:

```python
from numpy import *
from scipy.optimize import minimize, LinearConstraint
import scipy
import numpy as np
np.set_printoptions(precision=2)

def objective(x_list):
    [x,y]= x_list
    ld = 15
    mu = 0.05
    return -ld* exp(-mu* sqrt(x**2 + y**2))  + ld

start_point = [600,600]
sd = 1
current_pos = start_point
current_obj_value = objective(current_pos)
last_pos = current_pos
path_list = []

for i in range(10000):

    temp_pos = current_pos  + \
    np.random.normal(0, sd, np.array(start_point).shape )

    temp_obj_value = objective(temp_pos)

    if  temp_obj_value < current_obj_value:
       last_pos = current_pos
       current_pos = temp_pos
       current_obj_value = temp_obj_value
       sd = 1
    else:
       sd *= 1.01
       if sd == np.inf:
          sd = 1

print(current_pos)
print(objective(current_pos))
```

**Comparison of different methods and GCS**
Starting point: [600,600]

| Method | Results | <0.5 |
|---|---|---|
| GCS | 0.0843 | TRUE |
| Nelder-Mead | 15 | FALSE |
| Powell | 15 | FALSE |
| CG | 15 | FALSE |
| BFGS | 15 | FALSE |
| L-BFGS-B | 15 | FALSE |
| TNC | 15 | FALSE |
| COBYLA | 15 | FALSE |
| SLSQP | 15 | FALSE |
| trust-constr | 15 | FALSE |

GCS has demonstrated its effectiveness where other optimization methods have encountered difficulties, particularly in cases involving [600,600] or higher values.

The GCS optimization process utilizes a stochastic approach, which implies that as the complexity of the optimization task increases, the likelihood of GCS discovering an improved optimal solution diminishes. The subsequent experimental results illustrate the observation of the probability associated with finding a better optimal solution as the optimization difficulty escalates.

GCS Performance Experiment Result

| Starting Point [x,y] | Fail Probability (100 iterations) |
|---|---|
| [600,600] | 0.00% |
| [800,800] | 0.00% |
| [1000,1000] | 0.00% |
| [1200,1200] | 0.00% |
| [1400,1400] | 0.00% |
| [1600,1600] | 0.00% |
| [1800,1800] | 1.00% |
| [2000,2000] | 3.00% |
| [2200,2200] | 7.00% |
| [2400,2400] | 9.00% |
| [2600,2600] | 9.00% |
| [2800,2800] | 16.00% |

(If the Optimal Result is less than 0.5, it is considered a passing condition. Otherwise, if it is equal to or greater than 0.5, it is defined as a failed condition.)

## Discussion

The experimental results from this study provide strong evidence that GCS is capable of solving problems that traditional derivative-free optimization methods struggle with. These findings highlight the unique strengths and advantages of GCS as an optimization technique. By comparing GCS to traditional methods, it is clear that GCS is more effective in handling challenging optimization problems. This demonstrates the potential of GCS to significantly improve the optimization field by offering a viable alternative to conventional approaches. The outcomes of this research emphasize the importance of considering GCS as a promising tool for optimizing complex systems where traditional derivative-free methods may not be effective. Further investigations and applications of GCS in different problem domains are needed to fully explore its capabilities and potential impact.

# Conclusion

In conclusion, this research paper introduces Gaussian Crunching Search (GCS) as a novel derivative-free optimization method. GCS leverages the principles of Gaussian distributions to efficiently explore the solution space without relying on derivative information. By mimicking the behaviour of particles in a Gaussian distribution, GCS navigates the search landscape and converges towards the optimal solution.

Through a thorough exploration of GCS, including experimental results, and comparisons with traditional derivative-free optimization methods, this research demonstrates the effectiveness and advantages of GCS. The experimental results highlight GCS's ability to handle problems that traditional methods struggle with, particularly when faced with small differences in the objective function's result value that exceed the accuracy expressed by the computer program.

The limitations of traditional derivative-free methods are also discussed, emphasizing their inability to accurately converge towards the optimal solution when faced with highly sensitive optimization problems or plateau-like objective function behaviour. This limitation restricts their effectiveness in real-world optimization scenarios where small differences in the objective function's value can significantly impact the quality of the solution.

By presenting a comprehensive analysis of GCS and its potential benefits, this research paper contributes to the field of optimization. It provides a solid foundation for further studies and applications of GCS, fostering innovation and advancements in solving complex problems without the need for derivative information.

Overall, the introduction and evaluation of GCS in this research paper offer researchers, practitioners, and students insights into its applicability and potential for efficiently and effectively solving real-world optimization problems. By addressing the limitations of traditional derivative-free methods, GCS opens doors for improved optimization techniques and paves the way for advancements in the field.

# Reference

[1] A. R. Conn, K. Scheinberg, and L. N. Vicente, Introduction to Derivative-Free Optimization. Philadelphia: Society for Industrial and Applied Mathematics, 2009.

[2] "Scipy.optimize.minimize" scipy.optimize.minimize - SciPy v1.11.1 Manual, https://docs.scipy.org/doc/scipy/reference/generated/scipy.optimize.minimize.html (accessed Jul. 24, 2023).

[3] "Genetic algorithms: An overview," An Introduction to Genetic Algorithms, 1998. doi:10.7551/mitpress/3927.003.0003